\documentclass[final, reqno]{siamltex}

\usepackage{amsmath, amssymb, dsfont, mathrsfs}

\newtheorem{example}[theorem]{Example}
\newtheorem{remark}[theorem]{Remark}

\begin{document}

\title{Rate of Convergence of Weak Euler Approximation for Nondegenerate It\^{o} Diffusion and Jump Processes}

\author{Remigijus Mikulevi\v{c}ius \thanks{Department of Mathematics, University of Southern California, Los Angeles, USA \newline ({\tt mikulvcs@math.usc.edu}).}
        \and Changyong Zhang \thanks{Department of Mathematics, Uppsala University, Uppsala, Sweden \newline ({\tt changyong.zhang@math.uu.se}).}}

\maketitle

\begin{abstract}
The paper studies the rate of convergence of the weak Euler approximation for It\^{o} diffusion and jump processes with H\"{o}lder-continuous generators. It covers a number of stochastic processes including the nondegenerate diffusion processes and a class of stochastic differential equations driven by stable processes. To estimate the rate of convergence, the existence of a unique solution to the corresponding backward Kolmogorov equation in H\"{o}lder space is first proved. It then shows that the Euler scheme yields positive weak order of convergence.
\end{abstract}

\begin{keywords}
Stochastic differential equations, It\^{o} diffusion and jump processes, weak Euler approximation, rate of convergence, H\"{o}lder conditions, backward Kolmogorov equation
\end{keywords}

\section{Introduction}

\subsection{Continuous and Discontinuous Processes}

Since it was introduced in the early years of the twentieth century as a model for the physical phenomenon of Brownian motion by Einstein and Smoluchowski~\cite{Ein05, Smo06} and as a description of the dynamical evolution of stock prices by Bachelier~\cite{Bac00}, the Wiener process has been the most intensively-studied stochastic process.

Two important properties of the Wiener process are continuity of sample paths and scale invariance, while many phenomena that were first described by the Wiener process do not exhibit these properties. For example, the classical Black-Scholes model assumes that a stock price follows a geometric Brownian motion~\cite{BlS73}, whereas essentially stock prices change by units and move by jumps, further due to external shocks such as earnings announcements and natural disasters. Moreover, empirical studies have shown that stock prices follow distributions with heavy tails and skewness, which are incompatible with models based on the Wiener process.

On the other hand, stochastic processes with jumps, for instance L\'{e}vy processes~\cite{App04, Ber96, Sat99}, provide more flexibility in modeling large and sudden changes and naturally exhibit high skewness and leptokurtosis levels. This makes them useful tools in a broad variety of fields~\cite{BMR01}, for instance, differential geometry and extreme value theory in mathematics, Burgers' turbulence and quantum theory in physics, and particularly in finance~\cite{CoT12, Sch03}, from portfolio and risk management to option and bond pricing and hedging. Other applications can be found in engineering and science~\cite{CiP72, StK74, WoT05}.

\subsection{Analytic and Numerical Solutions}

Let $X = \{X_t\}_{t \in [0, T]}$, $T \in (0, \infty)$ be the process modeling the system of interest and $(\Omega, \mathcal{F}, P)$ be a complete probability space with a filtration $\mathbb{F} = \{\mathcal{F}_t\}_{t \in [0, T]}$ of $\sigma$-algebras satisfying the usual conditions. In many applications, it is frequently necessary to evaluate functionals of the underlying process, such as the moments. In particular, for a given test function $g$, the problem of computing $\mathrm{E}[g(X_T)]$ arises from various applications. For example, in telecommunication, the expectation represents the average energy of the system at time $T$, which is critical to the design and maintenance of telephone lines~\cite{StK74}. In finance and insurance, it is extremely common to evaluate $\mathrm{E}[g(X_T)]$, for instance, in pricing options and other derivatives~\cite{Kou02, Mer76}.

When both the model and test function are simple, as in the case that $g$ is sufficiently smooth and the increments of the driving process can be simulated, there exists a closed-form expression for the expectation. For instance in one dimension,
\begin{equation*}
X_t = X_0 (1 + c)^{N_t} \exp\big[\big(a - \frac{1}{2} b^2\big)t + b W_t\big]
\end{equation*}
is the solution to the stochastic differential equation
\begin{equation*}
X_t = X_0 + a \int_0^t X_{s-} ds + b \int_0^t X_{s-} dW_s + c \int_0^t X_{s - } dN_s, \ \forall t \in [0, T],
\end{equation*}
where $W = \{W_t\}_{t \in [0, T]}$ is a Wiener process, $N = \{N_t\}_{t \in [0, T]}$ is a Poisson process, and $X_{t-}$ denotes the left-hand limit of $X$ at $t$~\cite{Gla03}. In general, a system is more complicated and a closed-form solution is unrealistic.

An alternative option is to numerically approximate $\mathrm{E}[g(X_T)]$ by a discrete-time simulation of the stochastic process $X$~\cite{Gla03, TaR00}. This approach has been widely applied in practice~\cite{Ami93, CaH03, CoV05, Zha97}. One of the most commonly-used methods is Monte-Carlo simulation and the simplest scheme is the weak Euler approximation, for which it is of theoretical and practical importance to estimate the rate of convergence.

The Euler approximation $Y$ of a stochastic process $X$ is said to converge with a weak order $\kappa > 0$ if for any smooth function $g$, there exists a constant $K$, depending only on $g$, such that
\begin{equation*}
|\mathrm{E}[g(Y_T)] - \mathrm{E}[g(X_T)]| < K \delta^\kappa,
\end{equation*}
where $\delta > 0$ is the maximum step size of the time discretization.

\subsection{Smooth and H\"{o}lder-Continuous Coefficients}

The problem of estimating the rate of convergence has been systematically studied for processes with smooth coefficients. Milstein was one of the first who looked into the order of weak convergence of discrete-time approximations for diffusion processes~\cite{Mil78, Mil85}. Talay investigated a class of second weak-order approximations for diffusion processes~\cite{Tal84, Tal86}. Platen and Kloeden studied both Euler and higher-order schemes for It\^{o} processes~\cite{KlP92, Pla99}. For discrete-time approximations of It\^{o} processes with jumps, Mikulevi\v{c}ius and Platen showed first weak-order convergence in the case that the coefficient functions possess fourth-order continuous differentiability~\cite{MiP88}. Protter, Talay, Jacod, and Rubenthaler presented similar results for L\'{e}vy-driven stochastic differential equations~\cite{Jac04, PrT97, Rub03}.

In general, the coefficient and test functions do not necessarily satisfy the smoothness conditions assumed in the aforementioned papers. Mikulevi\v{c}ius and Platen first proved that there is still some weak-order convergence of the Euler approximation for diffusion processes under H\"{o}lder conditions~\cite{MiP911}. Kubilius and Platen generalized the result to diffusion processes with a finite number of jumps in finite time intervals~\cite{KuP02}.

In this paper, the rate of convergence for It\^{o} diffusion and jump processes under $\beta$-H\"{o}lder conditions is derived, by using the solution to the corresponding backward Kolmogorov equation associated with the underlying stochastic process. In Section~\ref{sec:model_diffusion_jump}, the model being considered is introduced and the main result is stated. In Section~\ref{sec:proof_diffusion_jump}, the main theorem is proved, with the necessary and essential technical results.

\section{Model and Result} \label{sec:model_diffusion_jump}

\subsection{It\^{o} Diffusion and Jump Process}

Denote $\mathds{R}_0^d = \mathds{R}^d \setminus \{0\}$ and $\tilde{\chi}^\alpha(y) = \mathbf{1}_{\{\alpha \in [1, 2]\}} \mathbf{1}_{\{|y| \le 1\}}$. For a fixed $\alpha \in (0, 2]$, consider a $d$-dimensional $\mathbb{F}$-adapted stochastic process $X = \{X_t\}_{t \in [0, T]}$, which solves
\begin{eqnarray} \label{eqn:diffusion_jump}
\begin{array}{rcl} X_t & = & {\displaystyle X_0 + \int_0^t a^\alpha(X_{s-})ds + \int_0^t b^\alpha(X_{s-})dW_s} \\
 & & + {\displaystyle \int_0^t \int_{\mathds{R}_0^d} y \big[\big(1 - \tilde{\chi}^\alpha(y)\big) p^\alpha(dy, ds) + \tilde{\chi}^\alpha(y) q^\alpha(dy, ds)\big]}, \end{array}
\end{eqnarray}
where $a^\alpha(x) = a(x) = \big(a^i(x)\big)_{1 \le i \le d}$ and $b^\alpha(x) = b(x) = \big(b^{ij}(x)\big)_{1 \le i, j \le d}$, $x \in \mathds{R}^d$ are measurable and bounded, $a = 0$ if $\alpha \in (0, 1)$ and $b = 0$ if $\alpha \ne 2$, and $b$ is nondegenerate; $W = \{W_t\}_{t \in [0, T]}$ is a $d$-dimensional $\mathbb{F}$-adapted standard Wiener process; $p^\alpha(dy, dt)$ is a jump measure on $\mathds{R}_0^d \times [0, T]$ with $p^\alpha(\Gamma, [0, t]) = \sum_{s \le t} \mathbf{1}_\Gamma(X_s - X_{s-})$, $\Gamma \in \mathcal{B}(\mathds{R}_0^d)$ and $q^\alpha(dy, dt) = p^\alpha(dy, dt) - \pi^\alpha(X_t, dy)dt$ is the corresponding martingale measure with $\pi^\alpha(x, dy)$, $x \in \mathds{R}^d$ being a measurable family of nonnegative measures on $\mathds{R}_0^d$.

For $\alpha \in (0, 2)$, it is assumed that $\pi^\alpha$ consists of two parts. The principal part has a nondegenerate density with respect to the L\'{e}vy measure of the spherically-symmetric $\alpha$-stable process $S = \{S_t\}_{t \in [0, T]}$ defined by
\begin{equation} \label{eqn:stable}
S_t = \int_0^t \int y \big[\big(1 - \chi^\alpha(y)\big) p_0(dy, ds) + \chi^\alpha(y) q_0(dy, ds)\big], \alpha \in (0, 2),
\end{equation}
where $\chi^\alpha(y) = \mathbf{1}_{\{\alpha = 1\}} \mathbf{1}_{\{|y| \le 1\}} + \mathbf{1}_{\{\alpha \in (1, 2)\}}$, $p_0(dy, ds)$ is the jump measure, and $\displaystyle q_0(dy, ds) = p_0(dy, ds) - \frac{dy}{|y|^{d + \alpha}}ds$ is the martingale measure. $S$ is a Wiener process if $\alpha = 2$. The second part has a density with respect to a lower-order L\'{e}vy measure. That is,
\begin{equation} \label{eqn:Levy_measure_diffusion_jump}
\pi^\alpha(x, dy) = m^\alpha(x, y) \frac{dy}{|y|^{d + \alpha}} + \rho^\alpha(x, y) \nu^\alpha(dy),
\end{equation}
where $\nu^\alpha$ is a nonnegative measure on $\mathds{R}_0^d$ and $m^\alpha, \rho^\alpha$ are nonnegative measurable functions such that $m^\alpha$ and $\int_{\mathds{R}_0^d} (|y|^\alpha \wedge 1) \rho^\alpha(x, y) \nu^\alpha(dy)$, $x \in \mathds{R}^d$ are bounded. For $\alpha = 2$, $\pi^\alpha(x, dy) = \rho^\alpha(x, y) \nu^\alpha(dy)$. In particular, if $\pi^\alpha = 0$, then $X$ is a diffusion process.

A large class of strong Markov processes satisfying (\ref{eqn:diffusion_jump}) has been constructed~\cite{AbK09, Kom84, LeM76, MiP923, Str75}. The processes are characterized by the drift coefficient $a$, diffusion coefficient $b$, and L\'{e}vy jump measure $\pi^\alpha(x, dy)$, or equivalently by their generators, as defined in Section~\ref{sec:proof_diffusion_jump}. They naturally arise in stochastic differential equations driven by L\'{e}vy processes. One example is the solution to the following stochastic differential equation driven by a spherically-symmetric stable process.

\begin{example} \label{exm:stable}
Let $S = \{S_t\}_{t \in [0, T]}$ be a $d$-dimensional standard spherically-symmetric $\alpha$-stable process defined by $(\ref{eqn:stable})$ and $L = (L^1, \ldots, L^d)$, with $L^i, i = 1, \dots, d$ being independent one-dimensional standard symmetric $\alpha^i$-stable processes independent of $S$. Consider for $t \in [0, T]$,
\begin{equation} \label{eqn:exm_stable}
X_t = X_0 + \mathbf{1}_{\{\alpha \in (1, 2]\}} \int_0^t a^\alpha(X_{s-})ds + \int_0^t c(X_{s-})dS_s + \int_0^t \mathrm{diag} \big(l(X_{s-})\big)dL_s,
\end{equation}
where $a^\alpha(x) = \big(a^i(x)\big)_{1 \le i \le d}$, $c(x) = \big(c^{ij}(x)\big)_{1 \le i, j \le d}$, and $l(x) = \big(l^i(x)\big)_{1 \le i \le d}$, $x \in \mathds{R}^d$ are $\beta$-H\"{o}lder continuous and bounded with $\beta > 0, \beta \ne \alpha, \beta \notin \mathds{N}$, and $\mathrm{diag}(l)$ is a $d \times d$-dimensional diagonal matrix with diagonal elements $l^1, \ldots, l^d$. It is assumed that $c$ is nondegenerate with $\inf_x \det{|c(x)|} > 0$ and $\alpha^i < \alpha \le 2$. In this case, $(\ref{eqn:Levy_measure_diffusion_jump})$ holds for $\pi^\alpha$ with
\begin{eqnarray} \label{eqn:exm_stable_measure}
\begin{array}{rcl} m^\alpha(x, y) & = & {\displaystyle \frac{|y|^{d + \alpha}}{|\det c(x)| |c(x)^{-1} y|^{d + \alpha}}}, \\
\rho^\alpha(x, y) & = & {\displaystyle \sum_{i = 1}^d \mathbf{1}_{\{y = y_i e_i\}} |l^i(x)|^{\alpha^i}}, \\
\nu^\alpha(dy) & = & {\displaystyle \sum_{i = 1}^d \mathbf{1}_{\{y = y_i e_i\}} \frac{dy_i}{|y_i|^{1 + \alpha^i}}}, \end{array}
\end{eqnarray}
where $\{e_i, i = 1, \ldots, d\}$ is the canonical basis of $\mathds{R}^d$.
\end{example}

Mathematical models defined by (\ref{eqn:diffusion_jump}) are used in fields such as finance and insurance to capture continuous and discontinuous uncertainty associated with various random dynamic phenomena.

\subsection{Weak Euler Approximation}

Let the time discretization $\{\tau\}_\delta = \{\tau_i\}_{i = 0, \dots, n_T}$ with a maximum step size $\delta \in (0, 1)$ be a partition of $[0, T]$ such that $0 = \tau_0 < \tau_1 < \dots < \tau_{n_T} = T$ and $\max_i (\tau_i - \tau_{i - 1}) \le \delta$. The weak Euler approximation of $X$ is an $\mathbb{F}$-adapted stochastic process $Y = \{Y_t\}_{t \in (0, T]}$ defined by
\begin{eqnarray} \label{eqn:Euler_diffusion_jump}
\begin{array}{rcl} Y_t & = & {\displaystyle X_0 + \int_0^t a^\alpha(Y_{\tau_{i_s}})ds + \int_0^t b^\alpha(Y_{\tau_{i_s}})dW_s} \\
 & & + {\displaystyle \int_0^t \int_{\mathds{R}_0^d} y \big[\big(1 - \tilde{\chi}^\alpha(y)\big) \tilde{p}^\alpha(dy, ds) + \tilde{\chi}^\alpha(y) \tilde{q}^\alpha(dy, ds)\big]}, \end{array}
\end{eqnarray}
where $\tau_{i_s} = \tau_i$ if $s \in [\tau_i, \tau_{i + 1})$, $i = 0, \ldots, n_T - 1$, $\tilde{p}^\alpha(dy, dt)$ is the jump measure and $\tilde{q}^\alpha(dy, dt) = \tilde{p}^\alpha(dy, dt) - \pi^\alpha(Y_{\tau_{i_t}}, dy)dt$ is the corresponding $\mathbb{F}$-adapted martingale measure on $\mathds{R}_0^d \times [0, T]$.

Contrary to those in (\ref{eqn:diffusion_jump}), the coefficients in (\ref{eqn:Euler_diffusion_jump}) are piecewise constant in each time interval $[\tau_i, \tau_{i + 1})$.

\begin{proposition}
For each stochastic process defined by $(\ref{eqn:diffusion_jump})$, there exists an Euler approximation defined by $(\ref{eqn:Euler_diffusion_jump})$.
\end{proposition}

\begin{proof}
By Lemma 14.50 in \cite{Jac79}, there exists a measurable function $l^\alpha: \mathds{R}^d \times \mathds{R}_0 \mapsto \mathds{R}^d$ such that $\displaystyle \pi^\alpha(x, dy) = \int_{\mathds{R}_0} \mathbf{1}_{dy} \big(l^\alpha(x, z)\big) \frac{dz}{z^2}, x \in \mathds{R}^d$. Given $X_t$ satisfying (\ref{eqn:diffusion_jump}), let $p^\prime(dz, dt)$ be an independent Poisson point measure on $\mathds{R}_0 \times [0, T]$ with a compensator $\displaystyle \frac{dz}{z^2} dt$, the weak Euler approximation for $t \in [0, T]$ is then defined by
\begin{eqnarray*}
\begin{array}{rcl} Y_t & = & {\displaystyle X_0 + \int_0^t a^\alpha(Y_{\tau_{i_s}})ds + \int_0^t b^\alpha(Y_{\tau_{i_s}})dW_s} \\
 & & + {\displaystyle \int_0^t \int_{\mathds{R}_0^d} l^\alpha(Y_{\tau_{i_s}}, z) \big[\big(1 - \tilde{\chi}^\alpha\big(l^\alpha(Y_{\tau_{i_s}}, z)\big)\big) p^\prime(dz, ds) + \tilde{\chi}^\alpha\big(l^\alpha(Y_{\tau_{i_s}}, z)\big) \big(p^\prime(dz, ds) - \frac{dz}{z^2} ds\big)\big]}, \end{array}
\end{eqnarray*}
where $\tilde{\chi}^\alpha(l^\alpha(Y_{\tau_{i_s}}, z)) = \mathbf{1}_{\{\alpha \in [1, 2]\}} \mathbf{1}_{\{|l^\alpha(Y_{\tau_{i_s}}, z)| \le 1\}}$.
\end{proof}

The result on rate of convergence is stated in Theorem~\ref{thm:main_diffusion_jump}, followed by the main assumptions \textup{A1} and \textup{A2}. If without being explicitly specified, $C = C(\cdot, \ldots, \cdot)$ denotes possibly different constants depending only on the corresponding arguments and the following notations are used.

Denote $\mathds{H} = [0, T] \times \mathds{R}^d$ and $\mathds{N} = \{0, 1, 2, \ldots\}$. For $x, y \in \mathds{R}^d$, write $(x, y) = \sum_{i = 1}^d x_i y_i$, $|x| = \sqrt{(x, x)}$ and $|B| = \sum_{i = 1}^d |B^{ii}|, B \in \mathds{R}^{d\times d}$.

For $(t, x) \in \mathds{H}$, multiindex $\gamma \in \mathds{N}^d$, and $i, j = 1, \dots, d$, denote
\begin{eqnarray*}
\partial_t u(t, x) & = & \frac{\partial}{\partial t} u(t, x), \quad \partial_x^\gamma u(t, x) = \frac{\partial^{|\gamma|}}{\partial^{\gamma_1} x_1 \dots \partial^{\gamma_d} x_d} u(t, x), \\
\partial_i u(t, x) & = & \frac{\partial}{\partial x_i} u(t, x), \quad \partial_{ij}^2 u(t, x) = \frac{\partial^2}{\partial x_i x_j} u(t, x), \\
\partial_x u(t, x) & = & \nabla u(t, x) = \big(\partial_1 u(t, x), \dots, \partial_d u(t, x)\big), \\
\partial^2 u(t, x) & = & \Delta u(t, x) = \sum_{i = 1}^d \partial_{ii}^2 u(t, x), \\
\partial^\alpha v(x) & = & \mathscr{F}^{-1}[|\xi|^\alpha \mathscr{F} v(\xi)](x), \alpha \in (0, 2),
\end{eqnarray*}
where $\mathscr{F}$ is the Fourier transform with respect to $x \in \mathds{R}^d$ and $\mathscr{F}^{-1}$ is the inverse Fourier transform: $\displaystyle \mathscr{F} v(\xi) = \int_{\mathds{R}^d} e^{-i(\xi, x)} u(x)dx$ and $\displaystyle \mathscr{F}^{-1} v(x) = \frac{1}{(2\pi)^d} \int_{\mathds{R}^d} e^{i(\xi, x)} v(\xi)d\xi$.

$C_{b}^\infty(\mathds{H})$ is the set of functions $u$ on $\mathds{H}$ such that $u(t, x), \forall t \in [0, T]$ is infinitely differentiable in $x$ and $\sup_{(t, x) \in \mathds{H}} |\partial_x^\gamma u(t, x)| < \infty, \forall \gamma \in \mathds{N}^d$. $C_0^\infty(G)$ is the set of infinitely differentiable functions on an open set $G \subseteq \mathds{R}^d$ with compact support.

For $\beta = [\beta]^ - + \{\beta\}^ + > 0$, where $[\beta]^ - \in \mathds{N}$ and $\{\beta\}^ + \in (0, 1]$, $C^\beta(\mathds{H})$ denotes the space of measurable functions $u$ on $\mathds{H}$ such that the norm
\begin{eqnarray*}
|u|_\beta & = & \sum_{|\gamma| \le [\beta]^-} \sup_{t, x} |\partial_x^\gamma u(t, x)| + \mathbf{1}_{\{\{\beta\}^ + < 1\}} \sup_{\substack{|\gamma| = [\beta]^ - , \\ t, x, h \ne 0}} \frac{|\partial_x^\gamma u(t, x + h) - \partial_x^\gamma u(t, x)|}{|h|^{\{\beta\}^ + }} \\
 & & + \mathbf{1}_{\{\{\beta\}^ + = 1\}} \sup_{\substack{|\gamma| = [\beta]^ - , \\ t, x, h \ne 0}} \frac{|\partial_x^\gamma u(t, x + h) - 2 \partial_x^\gamma u(t, x) + \partial_x^\gamma u(t, x - h)|}{|h|^{\{\beta\}^ + }}
\end{eqnarray*}
is finite. Accordingly, $C^\beta(\mathds{R}^d)$ denotes the corresponding space of functions on $\mathds{R}^d$. The classes $C^\beta$ are H\"{o}lder-Zygmund spaces and coincide with H\"{o}lder spaces if $\beta \notin \mathds{N}$~\cite{Tri92}.

Assume that $m^\alpha(x, y)$ and its partial derivatives $\partial_y^\gamma m^\alpha(x, y), |\gamma| \le d_0 = [\frac{d}{2}] + 1$ are continuous in $(x, y)$, $m^\alpha(x, y)$ is homogeneous in
$y$ with index zero, $m^\alpha \equiv 0$ for $\alpha = 2$, and
\begin{equation*} \label{eqn:sphere_integral}
\int_{S^{d - 1}} y m^\alpha(\cdot, y) \mu_{d - 1}(dy) = 0 \mbox{ for } \alpha = 1,
\end{equation*}
where $S^{d - 1}$ is the unit sphere in $\mathds{R}^d$ and $\mu_{d - 1}$ is the Lebesgue measure.

\begin{theorem} \label{thm:main_diffusion_jump}
Let $\alpha \in (0, 2]$, $\beta > 0, \beta \ne \alpha, \beta \notin \mathds{N}$, and $Y = \{Y_t\}_{t \in (0, T]}$ be the weak Euler approximation of the stochastic process $X = \{X_t\}_{t \in (0, T]}$ defined by $(\ref{eqn:diffusion_jump})$. Assume \textup{A1} and \textup{A2} hold. Then there exists a constant $C$ such that,
\begin{equation} \label{eqn:main_theorem_diffusion_jump}
|\mathrm{E}[g(Y_T)] - \mathrm{E}[g(X_T)]| \le C |g|_{\alpha + \beta} \delta^{\kappa(\alpha, \beta)}, \forall g \in C^{\alpha + \beta}(\mathds{R}^d),
\end{equation}
where
\begin{equation*}
\kappa(\alpha, \beta) = \left\{\begin{array}{cl}
\frac{\beta}{\alpha}, & \beta < \alpha, \\
1, & \beta > \alpha.
\end{array}
\right.
\end{equation*}
\end{theorem}

\begin{enumerate}
\item[\textbf{A1}] There exists a constant $\mu > 0$ such that for all $x \in \mathds{R}^d$ and $|\xi| = 1$,
 \begin{eqnarray} \label{eqn:assumption_bound_diffusion_jump}
 \begin{array}{rl} {\displaystyle \int_{S^{d - 1}} |(y, \xi)|^\alpha m^\alpha(x, y) dy \ge \mu}, & \alpha \in (0, 2), \\
 \big(B(x) \xi, \xi\big) \ge \mu, & \alpha = 2, \end{array}
 \end{eqnarray}
 where $B(x) = b(x)^* b(x), x \in \mathds{R}^d$ and
 \begin{equation} \label{eqn:assumption_bound_below_diffusion_jump}
 \lim_{\delta \downarrow 0} \sup_{x \in \mathds{R}^d} \int_{|y| \le \delta} |y|^\alpha \rho^\alpha(x, y) \nu^\alpha(dy) = 0;
 \end{equation}
\item[\textbf{A2}] For $\beta = [\beta] + \{\beta\} > 0$ with $[\beta] \in \mathds{N}$ and $\{\beta\} \in (0, 1)$,
 \begin{equation*}
 M_\beta^\alpha + N_\beta^\alpha < \infty,
 \end{equation*}
 where
 \begin{equation} \label{eqn:assumption_M_diffusion_jump}
 M_\beta^\alpha = \mathbf{1}_{\{\alpha = 1\}} |a|_\beta + \mathbf{1}_{\{\alpha = 2\}} |B|_\beta + \mathbf{1}_{\{\alpha \in (0, 2)\}} \sup_{\substack{|\gamma| \le d_0, \\ |y| = 1}} |\partial_y^\gamma m^\alpha(\cdot, y)|_\beta
 \end{equation}
 and
 \begin{eqnarray} \label{eqn:assumption_N_diffusion_jump}
 \begin{array}{rcl} N_\beta^\alpha & = & {\displaystyle \mathbf{1}_{\{\alpha \in (1, 2]\}} |a|_\beta \notag + \sup_{|\gamma| = [\beta], x} \int_{\mathds{R}_0^d} \big(|y|^\alpha \wedge 1\big) \big[\big|\rho^\alpha(x, y) \big| + \big|\partial_x^\gamma \rho^\alpha(x, y) \big|\big] \nu^\alpha(dy)} \\
 & & + {\displaystyle \sup_{\substack{|\gamma| = [\beta], \\ h \ne 0}} \frac{1}{|h|^{\beta - [\beta]}} \int_{\mathds{R}_0^d} \big(|y|^\alpha \wedge 1\big) \big|\partial_x^\gamma \rho^\alpha(x + h, y) - \partial_x^\gamma \rho^\alpha(x, y) \big|\nu^\alpha(dy)}. \end{array}
 \end{eqnarray}
\end{enumerate}

\begin{remark}
Assumptions \textup{A1} and \textup{A2} guarantee that the solution to the backward Kolmogorov equation associated with $X_t$ is $(\alpha + \beta)$-H\"{o}lder. They are in direct correspondence to the standard classical assumptions when the operator is differential. The regularity of the solution determines the rate of convergence of the weak Euler approximation.
\end{remark}

\begin{remark}
Since $g(X_T)$ is the value of the solution to the backward Kolmogorov equation $(\ref{eqn:Kolmogorov_back_diffusion_jump})$, a parabolic integro-differential equation of order $\alpha \in (0, 2]$, the estimate provides the rate of convergence of a probabilistic approximation $\mathrm{E}[g(Y_T)]$ to $\mathrm{E}[g(X_T)]$.
\end{remark}

\begin{remark}
Condition $(\ref{eqn:assumption_bound_diffusion_jump})$ holds with a constant $\mu > 0$ if, for example, there is a Borel set $\Gamma \subseteq S^{d - 1}$ such that $\mu_{d - 1}(\Gamma) > 0$ and $\displaystyle \inf_{x \in \mathds{R}^d, y \in \Gamma} m^\alpha(x, y) > 0$.

Condition $(\ref{eqn:assumption_bound_below_diffusion_jump})$ is satisfied if there is a measurable function $\rho^\alpha(y)$ such that $\rho^\alpha(x, y) \le \rho^\alpha(y)$ and $\displaystyle \int_{\mathds{R}_0^d}(|y|^\alpha \wedge 1) \rho^\alpha(y) \nu^\alpha(dy) < \infty$.

Assumption \textup{A2} holds for $\beta \notin \mathds{N}$ if and only if
\begin{equation*}
\mathbf{1}_{\{\alpha \in [1, 2]\}} |a|_\beta + \mathbf{1}_{\{\alpha = 2\}} |B|_\beta < \infty,
\end{equation*}
and there exists a constant $C$ such that for all multiindices $|\gamma| \le [\beta]$, $|\gamma^\prime| \le [\frac{d}{2}] + 1$ and $x \in \mathds{R}^d$, $\displaystyle \int_{\mathds{R}_0^d} \big(|y|^\alpha \wedge 1\big) \big|\partial_x^\gamma \rho^\alpha(x, y) \big| \nu^\alpha(dy) + |\partial_x^\gamma \partial_y^{\gamma^\prime} m^\alpha(x, y)| \le C$, and for all $x, \tilde{x} \in \mathds{R}^d, y \in S^{d - 1}$, and multiindices $|\gamma| = [\beta], |\gamma^\prime| \le [\frac{d}{2}] + 1$,
\begin{eqnarray*}
|\partial_x^\gamma \partial_y^{\gamma^\prime} m^\alpha(x, y) - \partial_x^\gamma \partial_y^{\gamma^\prime} m^\alpha(\tilde{x}, y)| & \le & C |x - \tilde{x} |^{\beta - [\beta]}, \\
\int_{\mathds{R}_0^d} \big(|y|^\alpha \wedge 1\big) \big|\partial_x^\gamma \rho^\alpha(x, y) - \partial_x^\gamma \rho^\alpha(\tilde{x}, y) \big| \nu^\alpha(dy) & \le & C |x - \tilde{x} |^{\beta - [\beta]}.
\end{eqnarray*}
The very last inequality holds if, for example, for all $x, \tilde{x} \in \mathds{R}^d, y \in \mathds{R}_0^d$,
\begin{equation*}
\big|\partial_x^\gamma \rho^\alpha(x, y) - \partial_x^\gamma \rho^\alpha(\tilde{x}, y) \big| \le C |x - \tilde{x} |^{\beta - [\beta]} \quad \mbox{and} \quad \int_{\mathds{R}_0^d}(|y|^\alpha \wedge 1) \nu^\alpha(dy) < \infty.
\end{equation*}
\end{remark}

\begin{remark}
For the process defined by $(\ref{eqn:exm_stable})$, assumptions \textup{A1} and \textup{A2} are satisfied if $c$ is nondegenerate with $\inf_x \det{|c(x)|} > 0$ and $c^{ij} \in C^\beta(\mathds{R}^d)$, $|l^i|^{\alpha^i} \in C^\beta(\mathds{R}^d)$, $i, j = 1, \ldots, d$.
\end{remark}

Applying Theorem~\ref{thm:main_diffusion_jump} to (\ref{eqn:exm_stable}) results in Corollary~\ref{cor:exm_stable_diffusion_jump}.

\begin{corollary} \label{cor:exm_stable_diffusion_jump}
Let $X = \{X_t\}_{t \in [0, T]}$ satisfy \textup{(\ref{eqn:exm_stable})} and $Y = \{Y_t\}_{t \in [0, T]}$ defined by
\begin{equation*}
Y_t = X_0 + \mathbf{1}_{\{\alpha \in (1, 2]\}} \int_0^t a^\alpha(Y_{\tau_{i_s}})ds + \int_0^t c(Y_{\tau_{i_s}})dS_s + \int_0^t \mathrm{diag} \big(l(Y_{\tau_{i_s}})\big)dL_s
\end{equation*}
be the weak Euler approximation. Then \textup{(\ref{eqn:main_theorem_diffusion_jump})} holds.
\end{corollary}

\begin{proof}
By applying Theorem $14.80$ in \textup{\cite{Jac79}} and changing the variable of integration, for any $f \in C_0^\infty(\mathds{R}_0^d)$, the compensator of
\begin{equation*}
\int_0^t \int f(y) p^\alpha(dy, ds) = \int_0^t \int f\big(c(X_{s-})y\big) p_0(dy, ds) + \sum_{i = 1}^d \int_0^t \int f\big(l^i(X_{s-})y_i e_i\big) p^{L_i}(dy_i, ds)
\end{equation*}
is
\begin{equation*}
\int_0^t \int f\big(c(X_{s-})y\big) \frac{dy}{|y|^{d + \alpha}} ds + \sum_{i = 1}^d \int_0^t \int f\big(l^i(X_{s-})y_i e_i\big) \frac{dy_i}{|y_i|^{1 + \alpha_i}} ds
\end{equation*}
or
\begin{equation*}
\int_0^t \int f(y) m^\alpha(X_s, y) \frac{dy}{|y|^{d + \alpha}} ds + \int_0^t \int f(y) \rho^\alpha(X_s, y) \nu^\alpha(dy) ds.
\end{equation*}
Thus,
\begin{equation*}
q^\alpha(dy, dt) = p^\alpha(dy, dt) - m^\alpha(X_t, y) \frac{dy}{|y|^{d + \alpha}} - \rho^\alpha(X_t, y) \nu^\alpha(dy)
\end{equation*}
is a martingale measure with $m^\alpha$, $\rho^\alpha$, and $\nu^\alpha$ being as defined in (\ref{eqn:exm_stable_measure}). Clearly, (\ref{eqn:diffusion_jump}) holds for $\alpha \in (0, 1)$. Since $m^\alpha(x, y)$ and $\rho^\alpha(x, y)$ are symmetric in $y$, then
\begin{equation*}
X_t = X_0 + \int_0^t \int_{|y| > 1} y p^\alpha(dy, ds) + \int_0^t \int_{|y| \le 1} y q^\alpha(ds, dy), \alpha \in [1, 2).
\end{equation*}
Hence, the statement follows by Theorem~\ref{thm:main_diffusion_jump}.
\end{proof}

\begin{remark}
Under the assumption of Example~\ref{exm:stable} with $\alpha = 2$, it was derived that the convergence rate of a diffusion process is of the order $\frac{1}{3 - \beta} < \kappa(2, \beta) = \frac{\beta}{2}$ for $\beta \in (1, 2)$~\textup{\cite{MiP911}}. Corollary~\ref{cor:exm_stable_diffusion_jump} improves the rate of convergence.
\end{remark}

\section{Proof} \label{sec:proof_diffusion_jump}

The rate of convergence is estimated by solving the associated backward Kolmogorov equation, whose operators are defined as follows.

For $u \in C^{\alpha + \beta}(\mathds{H})$ and $z \in \mathds{R}^d$, denote
\begin{eqnarray} \label{eqn:operator_A_diffusion_jump}
\begin{array}{rcl} \mathcal{A}_z^\alpha u(t, x) & = & {\displaystyle \mathbf{1}_{\{\alpha = 1\}} \big(a(z), \nabla_x u(t, x)\big) + \frac{1}{2}\mathbf{1}_{\{\alpha = 2\}} \sum_{i, j = 1}^d B^{ij}(z) \partial_{ij}^2 u(t, x)} \\
 & & + {\displaystyle \int_{\mathds{R}_0^d} \big[u(t, x + y) - u(t, x) - \chi^\alpha(y) \big(\nabla_x u(t, x), y\big)\big] m^\alpha(z, y) \frac{dy}{|y|^{d + \alpha}}}, \\
\mathcal{A}^\alpha u(t, x) & = & \mathcal{A}_x^\alpha u(t, x) = \mathcal{A}_z^\alpha u(t, x)|_{z = x}, \end{array}
\end{eqnarray}
where $\chi^\alpha(y) = \mathbf{1}_{\{\alpha = 1\}} \mathbf{1}_{\{|y| \le 1\}} + \mathbf{1}_{\{\alpha \in (1, 2)\}}$, and
\begin{eqnarray} \label{eqn:operator_B_diffusion_jump}
\begin{array}{rcl} \mathcal{B}_z^\alpha u(t, x) & = & {\displaystyle \mathbf{1}_{\{\alpha \in (1, 2]\}} \big(a(z) + \int_{|y| > 1} y m^\alpha(z, y) \frac{dy}{|y|^{d + \alpha}}, \nabla_x u(t, x)\big)} \\
 & & + {\displaystyle \int_{\mathds{R}_0^d} \big[u(t, x + y) - u(t, x) - \mathbf{1}_{\{\alpha \in (1, 2]\}} \mathbf{1}_{\{|y| \le 1\}} \big(\nabla_x u(t, x), y\big)\big] \rho^\alpha(z, y) \nu^\alpha(dy)}, \\
\mathcal{B}^\alpha u(t, x) & = & \mathcal{B}_x^\alpha u(t, x) = \mathcal{B}_z^\alpha u(t, x)|_{z = x}. \end{array}
\end{eqnarray}

The operator $\mathcal{L}^\alpha = \mathcal{A}^\alpha + \mathcal{B}^\alpha$ is the generator of $X_t$ defined in \textup{(\ref{eqn:diffusion_jump})}. $\mathcal{A}^\alpha$ is the principal part and $\mathcal{B}^\alpha$ is the lower-order or subordinated part.

\begin{remark}
If $a = 0$, $(B^{ij}) = I$, and $m^\alpha = 1$, then $\mathcal{A}^\alpha$ is the generator of a standard spherically-symmetric $\alpha$-stable process as defined in \textup{(\ref{eqn:stable})}.
\end{remark}

\begin{remark} \label{rmk:martingale_diffusion_jump}
Under assumptions \textup{A1} and \textup{A2}, there exists a unique weak solution to \textup{(\ref{eqn:diffusion_jump})} and the stochastic process
\begin{equation*}
u(X_t) - \int_0^t (\mathcal{A}^\alpha + \mathcal{B}^\alpha) u(X_s)ds, \forall u \in C^{\alpha + \beta}(\mathds{R}^d), \forall \beta > 0
\end{equation*}
is a martingale~\textup{\cite{MiP923}}.
\end{remark}

If $v(t, x), (t, x) \in \mathds{H}$ satisfies the backward Kolmogorov equation
\begin{eqnarray} \label{eqn:Kolmogorov_back_diffusion_jump}
\begin{array}{rcl} \big(\partial_t + \mathcal{A}_x^\alpha + \mathcal{B}_x^\alpha\big) v(t, x) & = & 0, \quad 0 \le t \le T, \\
v(T, x) & = & g(x), \end{array}
\end{eqnarray}
then by It\^{o}'s formula,
\begin{equation*}
\mathrm{E}[g(Y_T)] - \mathrm{E}[g(X_T)] = \mathrm{E}[v(T, Y_T) - v(0, Y_0)] = \mathrm{E}\big[\int_0^T (\partial_t + \mathcal{L}_{Y_{\tau_{i_s}}}^\alpha) v(s, Y_s) ds\big].
\end{equation*}

The regularity of $v$ determines a one-step estimate, a key step in estimating the rate of convergence. For $\beta \in (0, 1)$, the results on Kolmogorov equations in H\"{o}lder classes have been proved~\cite{MiP922, MiP09}. They can be extended to the case $\beta > 1$ in a standard analytic way. Due to the lack of regularity, probabilistic techniques such as stochastic flows cannot be applied. Instead, Fourier multipliers are used to estimate precisely the principal part of the operator in H\"{o}lder spaces, as shown in Lemma~\ref{lem:operator_A_bound} and Corollary~\ref{cor:operator_A_bound}. This also provides an alternative way to approach problems of the same kind.

\subsection{Backward Kolmogorov Equation}

In H\"{o}lder-Zygmund spaces, consider the backward Kolmogorov equation associated with $X_t$:
\begin{eqnarray} \label{eqn:Kolmogorov_diffusion_jump}
\begin{array}{rcl} \big(\partial_t + \mathcal{A}_x^\alpha + \mathcal{B}_x^\alpha\big) v(t, x) & = & f(t, x), \\
v(T, x) & = & 0. \end{array}
\end{eqnarray}
The regularity of its solution is essential for the one-step estimate, which determines the rate of convergence.

For $\beta \in (0, 1)$ and $f \in C^\beta(\mathds{H})$, it has been shown that there exists a unique solution $v \in C^{\alpha + \beta}(\mathds{H})$ to (\ref{eqn:Kolmogorov_diffusion_jump}) and $|v|_{\alpha + \beta} \le C |f|_\beta$~\cite{MiP09}. To extend the result to the case $\beta > 1$, a standard induction method is applied, with steps including considering the differences defining the derivatives, interpreting them as solutions to (\ref{eqn:Kolmogorov_diffusion_jump}), using uniform estimates, and passing to the limit. The result is stated in Theorem~\ref{thm:solution_Kolmogorov_diffusion_jump}.

\begin{theorem} \label{thm:solution_Kolmogorov_diffusion_jump}
Let $\alpha \in (0, 2]$, $\beta > 0, \beta \ne \alpha, \beta \notin \mathds{N}$, and $f \in C^\beta(\mathds{H})$. Assume \textup{A1} and \textup{A2} hold. Then there exist a unique solution $v \in C^{\alpha + \beta}(\mathds{H})$ to \textup{(\ref{eqn:Kolmogorov_diffusion_jump})} and a constant $C$ independent of $f$ such that $|u|_{\alpha + \beta} \le C |f|_\beta$.
\end{theorem}

An immediate consequence of Theorem~\ref{thm:solution_Kolmogorov_diffusion_jump} is the following statement.

\begin{corollary} \label{cor:solution_Cauchy_diffusion_jump}
Let $\alpha \in (0, 2]$ and $\beta > 0, \beta \ne \alpha, \beta \notin \mathds{N}$. Assume \textup{A1} and \textup{A2} hold. Then for $f \in C^\beta(\mathds{H})$ and $g \in C^{\alpha + \beta}(\mathds{R}^d)$, there exist a unique solution $v \in C^{\alpha + \beta}(\mathds{H})$ to the Cauchy problem
\begin{eqnarray} \label{eqn:Cauchy_diffusion_jump}
\begin{array}{rcl} \big(\partial_t + \mathcal{A}_x^\alpha + \mathcal{B}_x^\alpha\big) v(t, x) & = & f(t, x), \\
v(T, x) & = & g(x) \end{array}
\end{eqnarray}
and a constant $C$ independent of $f$ and $g$ such that $|v|_{\alpha + \beta} \le C(|f|_\beta + |g|_{\alpha + \beta})$.
\end{corollary}

To prove Theorem~\ref{thm:solution_Kolmogorov_diffusion_jump} and Corollary~\ref{cor:solution_Cauchy_diffusion_jump}, the H\"{o}lder-norm estimates of $\mathcal{A}^\alpha f$ and $\mathcal{B}^\alpha f$ are first derived for $f \in C^{\alpha + \beta}(\mathds{R}^d)$, $\beta > 0$. An auxiliary lemma about uniform convergence of H\"{o}lder-continuous functions is proved as well.

\subsubsection{Estimate of Principal Operator}

Let $\alpha \in (0, 2]$ and $\beta > 0$. For $z \in \mathds{R}^d$ and $f \in C^{\alpha + \beta}(\mathds{R}^d)$, the principal part of the generator is
\begin{equation} \label{eqn:operator_A_Fourier}
\mathcal{A}_z^\alpha f(x) = \mathscr{F}^{-1}[\psi^\alpha(z, \xi) \mathscr{F} f(\xi)](x), x \in \mathds{R}^d,
\end{equation}
where
\begin{eqnarray} \label{eqn:operator_A_Fourier_kernel}
\begin{array}{rcl} \psi^\alpha(z, \xi) & = & {\displaystyle - K \int_{S^{d - 1}} |(y, \xi)|^\alpha \big[1 - i\big(\mathbf{1}_{\{\alpha \ne 1\}} \tan{\frac{\alpha \pi}{2}} \text{sgn}(y, \xi)} \\
 & & - {\displaystyle \frac{2}{\pi} \mathbf{1}_{\{\alpha = 1\}} \text{sgn}(y, \xi) \ln{|(y, \xi)|}\big)\big] m^\alpha(z, y) \mu_{d - 1}(dy)} \\
 & & - {\displaystyle i\mathbf{1}_{\{\alpha = 1\}} \big(a(z), \xi\big) - \frac{1}{2} \mathbf{1}_{\{\alpha = 2\}} \big(B(z) \xi, \xi\big), z, \xi \in \mathds{R}^d} \end{array}
\end{eqnarray}
and $K = K(\alpha)$ is a constant depending on $\alpha$~\cite{Kom84, MiP922}. By Theorem 2.3.1 in \cite{SaT94}, for a fixed $z \in \mathds{R}^d$, $\exp \big\{t\psi^\alpha(z, \xi) \big\}$ is the characteristic function of a stable process, whose generator $\mathcal{A}_z^\alpha$ has the property that for $\alpha \in (0, 2)$,
\begin{eqnarray} \label{eqn:generator_stable}
\begin{array}{rcl} {\displaystyle \int \big[f(x + y) - f(x) - \chi^\alpha(y) \big(\nabla f(x), y\big)\big] \frac{dy}{|y|^{d + \alpha}}} & = & C(\alpha, d) \partial^\alpha f \\
 & = & C(\alpha, d) \mathscr{F}^{-1}[|\xi|^\alpha \mathscr{F} f(\xi)], \end{array}
\end{eqnarray}
where $\chi^\alpha(y) = \mathbf{1}_{\{\alpha = 1\}} \mathbf{1}_{\{|y| \le 1\}} + \mathbf{1}_{\{\alpha \in (1, 2)\}}$ and $C(\alpha, d)$ is a constant.

For the derivation of the estimate related to $\mathcal{A}^\alpha$, Fourier multipliers in $C^{\alpha + \beta}(\mathds{R}^d)$ is used~\cite{Tri83}. For $\beta > 0$, the H\"{o}lder-Zygmund space $C^\beta(\mathds{R}^d)$ coincides with the Besov space $B_{\infty \infty}^\beta$ and the theory of multipliers in the latter can be applied, by considering the equivalent norms in $C^\beta(\mathds{R}^d)$.

Let $\phi \in C_0^\infty(\mathds{R}^d)$ be a nonnegative function such that $\mathrm{supp} \phi = \{\xi: \frac{1}{2} \le |\xi| \le 2\}$ and $\displaystyle \sum_{j = - \infty}^\infty \phi (2^{- j} \xi) = 1, \forall \xi \ne 0$. Define $\varphi_k \in \mathcal{S}(\mathds{R}^d)$, $k = 0, \pm 1, \dots$ by
\begin{equation} \label{eqn:varphi_definition}
\mathscr{F} \varphi_k = \phi (2^{-k} \xi)
\end{equation}
and $\psi \in \mathcal{S}(\mathds{R}^d)$ by
\begin{equation} \label{eqn:psi_definition}
\mathscr{F} \psi = 1 - \sum_{k \ge 1} \mathscr{F} \varphi_k(\xi),
\end{equation}
where $\mathcal{S}(\mathds{R}^d)$ is the Schwartz space of rapidly-decaying smooth functions on $\mathds{R}^d$~\cite{BeL76}.

\begin{lemma} \label{rem:equivalent_norms}
For $\alpha \in (0, 2]$, $\beta > 0$, and $\gamma = \alpha + \beta$, it holds that
\begin{itemize}
\item[\textup{(i)}] $\displaystyle |u|_\beta \sim \sup_x |\psi * u(x)| + \sup_{k \ge 1} 2^{\beta k} \sup_x |\varphi_k * u(x)|$;
\item[\textup{(ii)}] $\displaystyle |u|_{\alpha, \beta} = |u|_0 + |\partial^\alpha u|_\beta \sim |u - \partial^\alpha u|_\beta \sim |u|_\gamma$.
\end{itemize}
\end{lemma}

Lemma~\ref{rem:equivalent_norms} can be proved by first defining a family of operators $J^s: \mathcal{S}^\prime(\mathds{R}^d) \to \mathcal{S}^\prime(\mathds{R}^d)$,
\begin{equation*}
J^s u = \mathscr{F}^{-1} \big((1 + |\cdot|^2)^{\frac{s}{2}} \mathscr{F} u\big) \quad \mbox{and} \quad I^s u = \mathscr{F}^{-1} \big((1 + |\cdot|^s) \mathscr{F} u\big), s > 0,
\end{equation*}
where $\mathcal{S}^\prime(\mathds{R}^d)$ is the Schwartz space of generalized functions. By Theorem 2.3.8 in \cite{Tri83}, $J^s: C^{\beta + s}(\mathds{R}^d) \to C^\beta(\mathds{R}^d)$, $\beta > 0$ is an isomorphism. By Proposition 2 of Section V.3.1 and Lemma 2 of Section V.3.2 in \cite{Ste70},
\begin{equation*}
\mathscr{F}^{-1} \big[\frac{1 + |\xi|^s}{(1 + |\xi|^2)^{\frac{s}{2}}} \mathscr{F} u\big] \quad \mbox{and} \quad \mathscr{F}^{-1} \big[\frac{(1 + |\xi|^2)^{\frac{s}{2}}}{1 + |\xi|^s} \mathscr{F} u\big], \ \ s > 0
\end{equation*}
map $L_{\infty}(\mathds{R}^d)$ onto $L_{\infty}(\mathds{R}^d)$. This implies that $I^s: C^{\beta + s}(\mathds{R}^d) \to C^\beta(\mathds{R}^d)$, $s > 0$, $\beta > 0$ is an isomorphism as well. The results then follows.

For the estimates of H\"{o}lder differences, the following Lemma~\ref{lem:g2f} from $2.6.1$ in \cite{Tri83} will be called.

\begin{lemma} \label{lem:g2f}
Let $\beta > 0$, $h \in C^{d_0}(\mathds{R}^d)$ for $d_0 = \big[\frac{d}{2} \big] + 1$, and $K_0$ be a constant such that $|\partial^\gamma h(\xi)| \le K_0(1 + |\xi|)^{-|\gamma|}$ for any $\xi \in \mathds{R}^d$ and every multiindex $\gamma$ with $|\gamma| \le d_0$. Then there exists a constant $C$ such that
\begin{equation*}
|\mathscr{F}^{-1} (h \mathscr{F} f)|_\beta \le C K_0 |f|_\beta, \forall f \in C^\beta(\mathds{R}^d).
\end{equation*}
\end{lemma}

For a solution $u \in C^{\alpha + \beta}(\mathds{H})$ to (\ref{eqn:Kolmogorov_diffusion_jump}) with $\beta > 1$, it is necessary to consider $\mathcal{A}^\alpha$ whose coefficients are differentiated. This requires estimating $\mathcal{A}^\alpha$ with coefficients not satisfying A1. Let $\bar{a} = (\bar{a}^i)_{1 \le i \le d} = \bar{a}(x)$, $\bar{b} = (\bar{b}^{ij})_{1 \le i, j \le d} = \bar{b}(x)$, $\bar{m}^\alpha = \bar{m}^\alpha(x, y)$, and $\bar{\rho}^\alpha = \bar{\rho}^\alpha(x, y)$, $x \in \mathds{R}^d, y \in \mathds{R}_0^d$ be measurable functions such that $\bar{m}^\alpha(x, y)$ and its partial derivatives $\partial_y^\gamma \bar{m}^\alpha(x, y), |\gamma| \le \big[\frac{d}{2} \big] + 1$ are continuous in $(x, y)$, $\bar{m}^\alpha(x, y)$ is
homogeneous in $y$ with index zero, and
\begin{equation*}
\int_{S^{d - 1}} y \bar{m}^\alpha(\cdot, y) \mu_{d - 1}(dy) = 0 \mbox{ for } \alpha = 1.
\end{equation*}

Define $\bar{M}_\beta^\alpha, \beta > 0$, $\bar{\mathcal{A}}^\alpha$, and $\bar{\psi}^\alpha$ by (\ref{eqn:assumption_M_diffusion_jump}), (\ref{eqn:operator_A_diffusion_jump}), and (\ref{eqn:operator_A_Fourier_kernel}), respectively, with $a$, $b$, and $m^\alpha$ replaced by $\bar{a}$, $\bar{b}$, and $\bar{m}^\alpha$, respectively. Apparently, for $\bar{\mathcal{A}}^\alpha$, the equality (\ref{eqn:operator_A_Fourier}) holds with $\bar{\psi}^\alpha$.

Let $\bar{\phi}(z, \xi) = \bar{\psi}^\alpha(z, \xi)(1 + |\xi|^\alpha)^{-1}, z, \xi \in \mathds{R}^d$. By Remark 10 in \cite{MiP09}, for every multiindex $\gamma$ with $|\gamma| \le d_0 = \big[\frac{d}{2} \big] + 1$, it holds that,
\begin{equation} \label{eqn:ineq_coeff}
\big|\partial_\xi^\gamma \bar{\phi}(\cdot, \xi) \big|_\beta \le C \bar{M}_\beta^\alpha |\xi|^{-|\gamma|}, \ \big|\partial_\xi^\gamma \bar{\phi}(\cdot, \xi) \big| \le C \bar{M}^\alpha |\xi|^{-|\gamma|}, \forall \beta > 0, \forall \xi \in \mathds{R}^d,
\end{equation}
where
\begin{equation*}
\bar{M}^\alpha = \sup_x\big[\mathbf{1}_{\{\alpha = 1\}} |\bar{a}(x)| + \mathbf{1}_{\{\alpha = 2\}} |\bar{b}(x)| + \mathbf{1}_{\{\alpha \in (0, 2)\}} \sup_{|\gamma| \le d_0, |y| = 1} |\partial_y^\gamma \bar{m}^\alpha(x, y)|\big].
\end{equation*}

Denote $\displaystyle \bar{\Phi} f(z, x) = \mathscr{F}_\xi^{-1} \big[\bar{\phi}(z, \xi) \mathscr{F} f(\xi)\big](x)$ and $\displaystyle \tilde{\Phi} f(x) = \bar{\Phi} f(x, x), z, x \in \mathds{R}^d$. The estimates for $\mathcal{A}^\alpha$ are derived from Lemma~\ref{lem:operator_A_bound}.

\begin{lemma} \label{lem:operator_A_bound}
Let $\beta > 0$ and $\bar{\beta} \in (0, \beta]$. Assume $\bar{M}_\beta^\alpha < \infty$. Then there exists a constant $C$ such that for all $f \in C^{\beta}(\mathds{R}^d)$,
\begin{eqnarray*}
|\bar{\Phi} f(z, \cdot)|_\beta & \le & C \bar{M}^\alpha |f|_\beta, z \in \mathds{R}^d, \\
|\bar{\Phi} f(\cdot, x)|_\beta & \le & C \bar{M}_\beta^\alpha |f|_{\bar{\beta}}, x \in \mathds{R}^d, \\
|\tilde{\Phi} f|_\beta & \le & C \big(\bar{M}^\alpha |f|_\beta + \bar{M}_\beta^\alpha |f|_{\bar{\beta}}\big).
\end{eqnarray*}
\end{lemma}

\begin{proof}
Let $\zeta_1 \in C_0^\infty(\mathds{R}^d)$ and $\zeta_2 = 1 - \zeta_1$ with $\zeta_1 \in [0, 1]$ and $\zeta_1(x) = 1$ if $|x| \le 1$. Then $\displaystyle \bar{\Phi} f(z, x) = \bar{\Phi}_1 f(z, x) + \bar{\Phi}_2f(z, x)$, where $\displaystyle \bar{\Phi}_k f(z, x) = \mathscr{F}^{-1} \big[\bar{\phi}(z, \xi) \zeta_k(\xi) \mathscr{F} f(\xi)\big](x), k = 1, 2$.

Clearly, $\bar{\Phi}_k f(z, x) = \eta_k(z, x) * \tilde{f}$, with
\begin{equation*}
\tilde{f} = \mathscr{F}^{-1} \big[(1 + |\xi|^\alpha)^{-1} \mathscr{F} f\big] \quad \mbox{and} \quad \eta_k(z, x) = \mathscr{F}^{-1} \big[\bar{\psi}^\alpha(z, \xi) \zeta_k(\xi)\big](x).
\end{equation*}

Let $\tilde{u} = \mathscr{F}^{-1} \zeta_1$, then $\tilde{u} \in \mathcal{S}(\mathds{R}^d)$ and
\begin{eqnarray*}
\eta_1(z, x) & = & \mathscr{F}^{-1} \big[\bar{\psi}^\alpha(z, \xi) \zeta_1(\xi)\big](x) \\
 & = & \mathscr{F}^{-1} \big[\bar{\psi}^\alpha(z, \xi) \mathscr{F} \tilde{u}(\xi)\big](x) \\
 & = & \mathbf{1}_{\{\alpha = 1\}} \big(\bar{a}(z), \nabla \tilde{u}(x)\big) + \frac{1}{2} \mathbf{1}_{\{\alpha = 2\}} \sum_{i, j = 1}^d \bar{B}^{ij}(z) \partial_{ij}^2 \tilde{u}(x) \\
 & & + \int \big[\tilde{u}(x + y) - \tilde{u}(x) - \chi^\alpha(y) \big(\nabla_x \tilde{u}(x), y\big)\big] \bar{m}^\alpha(z, y) \frac{dy}{|y|^{d + \alpha}}.
\end{eqnarray*}
Thus,
\begin{equation*}
\int|\eta_1(\cdot, x)|_\beta dx \le C \bar{M}_\beta^\alpha \quad \mbox{and} \quad \int|\eta_1(z, x)|dx \le C \bar{M}^\alpha.
\end{equation*}

By Lemma~\ref{rem:equivalent_norms}, $|\tilde{f} |_{\alpha + \beta} \le C |f|_\beta < \infty$. Then for any $x, z \in \mathds{R}^d$, $\bar{\beta} \le \alpha + \beta$,
\begin{eqnarray*}
\big|\int \eta_1(\cdot, y) \tilde{f}(x - y)dy\big|_\beta & \le & |\tilde{f} |_{\infty} \int|\eta_1(\cdot, y)|_\beta dy \le C \bar{M}_\beta^\alpha |f|_{\bar{\beta}}, \\
\big|\int \eta_1(z, y) \tilde{f}(\cdot - y)dy\big|_\beta & \le & |\tilde{f} |_\beta \int|\eta_1(z, y)|dy \le C \bar{M}^\alpha |f|_\beta.
\end{eqnarray*}
Hence,
\begin{equation*}
|\bar{\Phi}_1 f(\cdot, x)|_\beta \le C \bar{M}_\beta^\alpha |f|_{\bar{\beta}} \quad \mbox{and} \quad |\bar{\Phi}_1 f(z, \cdot)|_\beta \le C \bar{M}^\alpha |f|_\beta.
\end{equation*}

Lemma~\ref{lem:g2f} and (\ref{eqn:ineq_coeff}) imply that for any $\bar{\beta} \in (0, \beta]$,
\begin{equation*}
|\bar{\Phi}_2 f(\cdot, x)|_\beta \le C \bar{M}_\beta^\alpha |f|_{\bar{\beta}} \quad \mbox{and} \quad |\bar{\Phi}_2 f(z, \cdot)|_\beta \le C \bar{M}^\alpha |f|_\beta.
\end{equation*}

For any multiindex $|\gamma| \le [\beta]^ - $,
\begin{equation*}
\partial^\gamma [\tilde{\Phi}_k f(x, x)] = \sum_{\mu + \nu = \gamma} \partial_x^\mu \partial_z^\nu \tilde{\Phi}_k f(z, x)|_{z = x} = \sum_{\mu + \nu = \gamma} \partial_z^\nu \tilde{\Phi}_k(\partial^\mu f)(z, x)|_{z = x}, k = 1, 2. \notag
\end{equation*}
It then follows that $\displaystyle|\tilde{\Phi}_k f|_\beta \le C(\bar{M}^\alpha |f|_\beta + \bar{M}_\beta^\alpha |f|_{\bar{\beta}}) |f|_\beta, k = 1, 2$, where $\tilde{\Phi}_k f(x) = \tilde{\Phi}_k f(x, x)$, $\forall x \in \mathds{R}^d$.
\end{proof}

\begin{corollary} \label{cor:operator_A_bound}
Let $\alpha \in (0, 2]$, $\beta > 0$, $\bar{\beta} \in (0, \beta]$, and $\bar{M}_\beta^\alpha < \infty$. Assume $f \in C^{\alpha + \beta}(\mathds{R}^d)$. Then there exists a constant $C$ such that
\begin{eqnarray*}
|\bar{\mathcal{A}}_\cdot^\alpha f(x)|_\beta & \le & C \bar{M}_\beta^\alpha |f|_{\alpha + \bar{\beta}}, x \in \mathds{R}^d, \\
|\bar{\mathcal{A}}_z^\alpha f(\cdot)|_\beta & \le & C \bar{M}^\alpha |f|_{\alpha + \beta}, z \in \mathds{R}^d, \\
|\bar{\mathcal{A}}^\alpha f|_\beta & \le & C (\bar{M}_\beta^\alpha |f|_{\alpha + \beta} + \bar{M}^\alpha |f|_{\alpha + \beta}).
\end{eqnarray*}
\end{corollary}

\begin{proof}
Let $\tilde{f} = \mathscr{F}^{-1}[(1 + |\xi|^\alpha) \mathscr{F} f]$. Then by Lemma~\ref{rem:equivalent_norms}, $\tilde{f} \in C^\beta(\mathds{R}^d)$, $|\tilde{f} |_\beta \le C |f|_{\alpha + \beta}$, and $\displaystyle \bar{\mathcal{A}}_z^\alpha f(x) = \mathscr{F}^{-1}[\bar{\psi}^\alpha(z, \xi)(1 + |\xi|^\alpha)^{-1} \mathscr{F} \tilde{f}] = \mathscr{F}^{-1}[\bar{\phi}(z, \xi) \mathscr{F} \tilde{f}]$. The statement follows by Lemma~\ref{lem:operator_A_bound}.
\end{proof}

\subsubsection{Estimate of Lower-Order Operator}

As in the case of $\mathcal{A}^\alpha$, it is necessary to consider $\mathcal{B}^\alpha$ whose coefficients are differentiated. Let $\bar{a} = (\bar{a}^i)_{1 \le i \le d} = \bar{a}(x)$ and $\bar{\rho}^\alpha = \bar{\rho}^\alpha(x, y)$, $x \in \mathds{R}^d$, $y \in \mathds{R}_0^d$ be measurable functions. Define $\bar{N}_\beta^\alpha$ and $\mathcal{\bar{B}}^\alpha$ by (\ref{eqn:assumption_N_diffusion_jump}) and (\ref{eqn:operator_B_diffusion_jump}), respectively, with $a$ and $\rho^\alpha$ being replaced by $\bar{a}$ and $\bar{\rho}^\alpha$.

Lemma~\ref{lem:difference_representation}, the result of which is stated in Lemma 2.1 in \cite{Kom84}, will be called for the estimate of $\mathcal{B}^\alpha$.

\begin{lemma} \label{lem:difference_representation}
Let $\delta \in (0, 1)$ and $u \in C_0^{\infty}(\mathds{R}^d)$. Denote
\begin{equation*}
k^{(\delta)}(y, z) = |z + y|^{-d + \delta} - |z|^{-d + \delta}.
\end{equation*}
Then there exist constants $C$ and $K = K(\delta, d)$ such that
\begin{equation*}
\int |k^{(\delta)}(y, z)| dz \le C |y|^\delta, \forall y \in \mathds{R}^d
\end{equation*}
and
\begin{equation*}
u(x + y) - u(x) = K \int k^{(\delta)}(y, z) \partial^\delta u(x - z) dz.
\end{equation*}
\end{lemma}

Denote $\displaystyle \bar{N}^\alpha = \mathbf{1}_{\{\alpha \in (1, 2]\}} \sup_x |\bar{a}(x)|_\beta + \sup_{x \in \mathds{R}^d} \int_{\mathds{R}_0^d} \big(|y|^\alpha \wedge 1\big) \big|\bar{\rho}^\alpha(x, y) \big| \nu^\alpha(dy)$.

\begin{lemma} \label{lem:operator_B_bound}
Let $\beta > 0$, $\bar{\beta} \in (0, \beta]$, and $\bar{N}_\beta^\alpha < \infty$. Then there exists a constant $C$ such that for all $f \in C^{\alpha + \beta}(\mathds{R}^d)$,
\begin{eqnarray*}
|\bar{\mathcal{B}}_z^\alpha f(\cdot)|_\beta & \le & C \bar{N}^\alpha |f|_{\alpha + \beta}, z \in \mathds{R}^d, \label{ff15} \\
|\bar{\mathcal{B}}_\cdot^\alpha f(x)|_\beta & \le & C \bar{N}_\beta^\alpha |f|_{\alpha + \bar{\beta}}, x \in \mathds{R}^d, \label{ff16} \\
|\bar{\mathcal{B}}^\alpha f|_\beta & \le & C \bar{N}_\beta^\alpha |f|_{\alpha + \beta}. \notag
\end{eqnarray*}
\end{lemma}

\begin{proof}
Rewrite $\bar{\mathcal{B}}_z^\alpha f(x) = \bar{\mathcal{B}}_z^{\alpha, 1} f(x) + \bar{\mathcal{B}}_z^{\alpha, 2} f(x)$, where
\begin{eqnarray*}
\bar{\mathcal{B}}_z^{\alpha, 1} f(x) & = & \int_{|y| > 1} [f(x + y) - f(x)] \bar{\rho}^\alpha(z, y) \nu^\alpha(dy) + \mathbf{1}_{\{\alpha \in (1, 2]\}} \big(\bar{a}(z), \nabla f(x)\big), \\
\bar{\mathcal{B}}_z^{\alpha, 2} f(x) & = & \mathbf{1}_{\{\alpha \in (0, 1]\}} \int_{|y| \le 1}[f(x + y) - f(x)] \bar{\rho}^\alpha(z, y) \nu^\alpha(dy) \\
 & & + \mathbf{1}_{\{\alpha \in (1, 2]\}} \int_{|y| \le 1} \int_0^1 \big(\nabla f(x + sy) - \nabla f(x), y\big) \bar{\rho}^\alpha(z, y)ds \nu^\alpha(dy).
\end{eqnarray*}
By Lemma~\ref{lem:difference_representation},
\begin{eqnarray*}
\bar{\mathcal{B}}_z^{\alpha, 2} f(x) & = & \mathbf{1}_{\{\alpha \in (0, 1)\}} K \int_{|y| \le 1} \int k^{(\alpha)}(y, \tilde{y}) \partial^\alpha f(x - \tilde{y}) d\tilde{y} \bar{\rho}^\alpha(z, y) \nu^\alpha(dy) \notag \\
 & & + \mathbf{1}_{\{\alpha = 1\}} \int_{|y| \le 1} \big(\nabla f(x + sy), y\big) \bar{\rho}^\alpha(z, y) \nu^\alpha(dy) \notag \\
 & & + \mathbf{1}_{\{\alpha \in (1, 2)\}} K \int_{|y| \le 1} \int_0^1 \big(\int k^{(\alpha - 1)}(sy, \tilde{y}) \partial^{\alpha - 1} \nabla f(x - \tilde{y}) d\tilde{y}, y\big) \times \bar{\rho}^\alpha(z, y) ds \nu^\alpha(dy) \notag \\
 & & + \mathbf{1}_{\{\alpha = 2\}} \sum_{i, j = 1}^d \int_{|y| \le 1} \int_0^1 (1 - s) \partial_{ij}^2 f(x + sy) y_i y_j ds \bar{\rho}^2(z, y) \nu^2(dy).
\end{eqnarray*}

For $\alpha > 1$, by Lemma~\ref{rem:equivalent_norms},
\begin{equation*}
|\partial^{\alpha - 1} \nabla f|_\beta \le |(1 - \partial^{\alpha - 1}) \nabla f|_\beta + |\nabla f|_\beta \le C (|\nabla f|_{\alpha + \beta - 1} + |\nabla f|_\beta) \le C |f|_{\alpha + \beta}, \forall \beta > 0.
\end{equation*}
The first two inequalities then follow. For example, if $\alpha \in (1, 2), \beta \notin \mathds{N}$,
\begin{eqnarray*}
 & & \big|\int_{|y| \le 1} \int_0^1 \big(\int k^{(\alpha - 1)}(sy, \tilde{y}) \partial^{\alpha - 1} \nabla f(x - \tilde{y})d\tilde{y}, y\big) \bar{\rho}^\alpha(\cdot, y) ds \nu^\alpha(dy) \big|_\beta \\
 & \le & C \sup_x |\partial^{\alpha - 1} \nabla f(x)| \sup_{\substack{|\gamma| \le [\beta], \\z, h \ne 0}} \int_{|y| \le 1} |y|^\alpha \big[|\partial^\gamma \bar{\rho}^\alpha(z, y)| + \frac{|\partial^\gamma \bar{\rho}^\alpha(z + h, y)| - |\partial^\gamma \bar{\rho}^\alpha(z, y)|}{|h|^{\beta - [\beta]}} \big] \nu^\alpha(dy),
\end{eqnarray*}
and
\begin{eqnarray*}
 & & \big|\int_{|y| \le 1} \int_0^1 \big(\int k^{(\alpha - 1)}(sy, \tilde{y}) \partial^{\alpha - 1} \nabla f(\cdot - \tilde{y})d\tilde{y}, y\big) \bar{\rho}^\alpha(z, y)ds \nu^\alpha(dy) \big|_\beta \\
 & \le & C |\partial^{\alpha - 1} \nabla f|_\beta \int_{|y| \le 1} |y|^\alpha |\bar{\rho}^\alpha(z, y)| \nu^\alpha(dy).
\end{eqnarray*}

Also, for $|\gamma| \le [\beta]^ - $, $\displaystyle \partial^\gamma (\bar{\mathcal{B}}^\alpha f) = \sum_{\kappa + \mu = \gamma} \partial_z^{\kappa} \bar{\mathcal{B}}_z^{\alpha, k}(\partial^\mu f)(x)|_{z = x}, k = 1, 2$. The third inequality follows as well.
\end{proof}

\subsubsection{Uniform Convergence of H\"{o}lder-Continuous Functions}

The result of Lemma~\ref{lem:P2Limit} will be applied in the proof of Theorem~\ref{thm:solution_Kolmogorov_diffusion_jump} by induction and passing to the limit.

\begin{lemma} \label{lem:P2Limit}
Assume $u_n \in C^\beta(\mathds{R}^d)$, $n \in \mathds{N}$ with $\sup_n |u_n|_\beta < \infty$ and $u_n \to u$ uniformly on compact subsets. Then $u \in C^\beta$, $|u|_\beta \le \sup_n|u_n|_\beta$, and $\displaystyle \partial^\delta \partial^\gamma u_n \to \partial^\delta \partial^\gamma u, |\gamma| \le [\beta]^ - $ uniformly on compact subsets as $n \to \infty$ for any $\delta \in [0, 1)$ such that $[\beta]^ - + \delta < \beta$.
\end{lemma}

\begin{proof}
Let $\psi$, $\varphi_k$ be functions defined by (\ref{eqn:varphi_definition}) and (\ref{eqn:psi_definition}), respectively. If $u_n \to u$ uniformly on compact sets, then
\begin{equation*}
|\psi * u(x)| = \lim_n |\psi * u_n(x)| \le \sup_n \sup_y |\psi * u_n(y)|, \forall x \in \mathds{R}^d
\end{equation*}
and
\begin{equation*}
2^{\beta k} |\varphi_k * u(x)| = 2^{\beta k} \lim_n |\varphi_k * u_n(x)| \le \sup_n \sup_k 2^{\beta k} \sup_y |\varphi_k * u_n(y)|, \forall x \in \mathds{R}^d.
\end{equation*}
Hence by Lemma~\ref{rem:equivalent_norms} (i), $|u|_\beta \le \sup \limits_n |u_n|_\beta < \infty$. By the Arzel\`{a}-Ascoli theorem, there exist continuous functions $v_\gamma(x), x \in \mathds{R}^d, |\gamma| \le [\beta]^ - $ and a subsequence $u_{n_k}$, whose limit is continuously differentiable up to $[\beta]^ - $, such that $\partial^\gamma u_{n_k} \to v_\gamma$ uniformly on compact subsets of $\mathds{R}^d$ as $k \to \infty$. Then by Theorem 3.6.1 in \cite{Car67}, $v_\gamma = \partial^\gamma v_0$. By (\ref{eqn:generator_stable}), for $\delta \in [0, 1)$ such that $[\beta]^ - + \delta < \beta$ and $|\mu| \le [\beta]^ - $, $\displaystyle \partial^\delta \partial^\mu u_n(x) = K \int \big[\partial^\mu u_n(x + y) - \partial^\mu u_n(x)\big] \frac{dy}{|y|^{d + \delta}}$. Passing to the limit yields that $\partial^\delta \partial^\mu u_n \to \partial^\delta \partial^\mu u$ uniformly on compact subsets as $n \to \infty$.
\end{proof}

\subsubsection{Solution to Kolmogorov Equation}

Theorem~\ref{thm:solution_Kolmogorov_diffusion_jump} is proved by induction.

\begin{proof}
For $\alpha \in (0, 2]$ and $\beta \in (0, 1)$, given $f \in C^\beta(\mathds{H})$, there exists a unique solution $u \in C^{\alpha + \beta}(\mathds{H})$ to the Kolmogorov equation (\ref{eqn:Kolmogorov_diffusion_jump}) and $|u|_{\alpha + \beta} \le C |f|_\beta$~\cite{MiP922}.

Assume the result holds for $\beta \in \bigcup\limits_{l = 0}^{n - 1}(l, l + 1)$, $n \in \mathds{N}$. Suppose that $\beta \in (n, n + 1)$ and $f \in C^\beta$. Let $\bar{\beta} = \beta - 1$. Then $\bar{\beta} \in (n - 1, n)$, $f \in C^{\bar{\beta}}(\mathds{H})$ as well, and there exists a unique solution $v \in C^{\alpha + \bar{\beta}}(\mathds{H})$, $\alpha \in (0, 2]$ to the Cauchy problem (\ref{eqn:Kolmogorov_diffusion_jump}) and $|v|_{\alpha + \bar{\beta}} \le C |f|_{\bar{\beta}}$.

For $h \in \mathds{R}$ and $k = 1, \dots, d$, denote
\begin{equation*}
v_k^h(t, x) = \frac{v(t, x + h e_k) - v(t, x)}{h},
\end{equation*}
where $\{e_k, k = 1, \ldots, d\}$ is the canonical basis in $\mathds{R}^d$. Let
\begin{equation*}
\mathcal{A}_z^{\alpha, k, h} v(t, x) = \frac{1}{h} \big(\mathcal{A}_{z + h e_k}^\alpha - \mathcal{A}_z^\alpha\big) v(t, x),
\end{equation*}
and
\begin{equation*}
\mathcal{B}_z^{\alpha, k, h} v(t, x) = \frac{1}{h} \big(\mathcal{B}_{z + h e_k}^\alpha - \mathcal{B}_z^\alpha\big) v(t, x).
\end{equation*}

Apparently,
\begin{eqnarray} \label{eqn:cauchy_1to2_h}
\begin{array}{rcl} {\displaystyle \big(\partial_t + \mathcal{A}_{x + h e_k}^\alpha + \mathcal{B}_{x + h e_k}^\alpha\big) v(t, x + h e_k)} & = & f(t, x + h e_k), \\
v(T, x + h e_k) & = & 0, \ \ k = 1, \dots, d. \end{array}
\end{eqnarray}

Subtracting (\ref{eqn:Kolmogorov_diffusion_jump}) from (\ref{eqn:cauchy_1to2_h}) and dividing the difference by $h$ yields
\begin{eqnarray} \label{eqn:cauchy_diff}
\begin{array}{rcl} {\displaystyle \big(\partial_t + \mathcal{A}_x^\alpha + \mathcal{B}_x^\alpha\big) v_k^h(t, x)} & = & {\displaystyle f_k^h(t, x) - \mathcal{A}_x^{\alpha, k, h} v(t, x + h e_k) - \mathcal{B}_x^{\alpha, k, h} v(t, x + h e_k)}, \\
v_k^h(T, x) & = & 0, \ \ k = 1, \dots, d. \end{array}
\end{eqnarray}

Since $f \in C^\beta(\mathds{H})$ and
\begin{equation*}
f_k^h(t, x) = \frac{f(t, x + h e_k) - f(t, x)}{h} = \int_0^1 \partial_k f(t, x + h e_k s)ds, \forall h \ne 0,
\end{equation*}
then
\begin{equation} \label{eqn:norm_f_bound}
|f_k^h|_{\bar{\beta}} \le C |\nabla f|_{\beta - 1} \le C |f|_\beta,
\end{equation}
where $C$ is a constant independent of $h$. Since $v \in C^{\alpha + \bar{\beta}}(\mathds{H})$, then $v_k^h \in C^{\alpha + \bar{\beta}}(\mathds{H})$. For $x \in \mathds{R}^d$, let
\begin{eqnarray*}
\bar{a}_{h, k}(x) & = & \frac{1}{h} \big(a(x + h e_k) - a(x)\big) = \int_0^1 \partial_k a(x + h e_k s)ds, \\
\bar{B}_{h, k}^{ij}(x) & = & \int_0^1 \partial_k B^{ij}(x + s h e_k)ds, \\
\bar{m}_{h, k}^\alpha(x, y) & = & \int_0^1 \partial_k m^\alpha(x + h e_k s, y)ds, \\
\bar{\rho}_{h, k}^\alpha(x, y) & = & \int_0^1 \partial_k \rho^\alpha(x + s h e_k, y)ds, y \in \mathds{R}_0^d.
\end{eqnarray*}
Then for $(t, x) \in \mathds{H}$,
\begin{eqnarray*}
 & & \mathcal{A}_x^{\alpha, k, h} v(t, x + h e_k) \\
 & = & \frac{1}{h} \big(\mathcal{A}_{x + h e_k}^\alpha - \mathcal{A}_x^\alpha\big) v(t, x + h e_k) \\
 & = & \mathbf{1}_{\{\alpha = 1\}} \big(\bar{a}_{h, k}(x), \nabla v(t, x + h e_k)\big) + \frac{1}{2} \mathbf{1}_{\{\alpha = 2\}} \sum_{i, j = 1}^d \bar{B}_{h, k}^{ij}(x) \partial_{ij}^2 v(t, x + h e_k) \\
 & & + \int \big[v(t, x + h e_k + y) - v(t, x + h e_k) - \chi^\alpha(y) \big(\nabla_x v(t, x + h e_k + y), y\big)\big] \bar{m}_{h, k}^\alpha(x, y) \frac{dy}{|y|^{d + \alpha}}
\end{eqnarray*}
and
\begin{eqnarray*}
\mathcal{B}_x^{\alpha, k, h} v(t, x + h e_k) & = & \frac{1}{h} \big(\mathcal{B}_{x + h e_k}^\alpha - \mathcal{B}_x^\alpha\big) v(t, x + h e_k) \\
 & = & \mathbf{1}_{\{\alpha \in (1, 2]\}} \big(\bar{a}_{h, k}(x), \nabla v(t, x + h e_k)\big) \\
 & & + \int_{|y| > 1} \nabla_y^1 v(t, x + h e_k) \bar{\rho}_{h, k}^\alpha(x, y) \nu^\alpha(dy) \\
 & & + \mathbf{1}_{\{\alpha \in (1, 2]\}} \int_{|y| \le 1} \nabla_y^2 v(t, x + h e_k) \bar{\rho}_{h, k}^\alpha(x, y) \nu^\alpha(dy) \\
 & & + \mathbf{1}_{\{\alpha \in (0, 1]\}} \int_{|y| \le 1} \nabla_y^1 v(t, x + h e_k) \bar{\rho}_{h, k}^\alpha(x, y) \nu^\alpha(dy),
\end{eqnarray*}
where
\begin{eqnarray*}
\nabla_y^1 v(t, x) & = & v(t, x + y) - v(t, x), \\
\nabla_y^2 v(t, x) & = & v(t, x + y) - v(t, x) - \big(\nabla v(t, x), y\big).
\end{eqnarray*}

By applying Corollary~\ref{cor:operator_A_bound} with $\bar{a} = \bar{a}_{h, k}$, $\bar{B} = \bar{B}_{h, k}$, $\bar{m}^\alpha = \bar{m}_{h, k}^\alpha$, $\beta = \bar{\beta}$ and by the induction assumption, it follows that
\begin{equation*}
|\mathcal{A}^{\alpha, k, h} v|_{\bar{\beta}} \le C M_\beta^\alpha |v|_{\alpha + \bar{\beta}} \le C M_\beta^\alpha |f|_{\bar{\beta}}, k = 1, \dots, d,
\end{equation*}
with a constant $C$ independent of $h$ and $f$.

Applying Lemma~\ref{lem:operator_B_bound} with $\bar{a} = \bar{a}_{h, k}$, $\bar{\rho}^\alpha = \bar{\rho}_{h, k}^\alpha$, $\beta = \bar{\beta}$, $f = v(t, \cdot + h e_k)$ together with the induction assumption yields
\begin{equation} \label{eqn:norm_B_bound}
\big|\mathcal{B}^{\alpha, k, h} v\big|_{\bar{\beta}} \le C N_\beta^\alpha |v|_{\alpha + \bar{\beta}} \le C N_\beta^\alpha |f|_{\bar{\beta}}, k = 1, \dots, d,
\end{equation}
where $C$ is a constant independent of $h$ and $f$.

Hence, $f_k^h(t, x) - \mathcal{A}_x^{\alpha, k, h} v(t, x + h e_k) - \mathcal{B}_x^{\alpha, k, h} v(t, x + h e_k) \in C^{\bar{\beta}}(\mathds{H})$ and $v_k^h \in C^{\alpha + \bar{\beta}}(\mathds{H})$ satisfies (\ref{eqn:cauchy_diff}). By the induction assumption and (\ref{eqn:norm_f_bound}) - (\ref{eqn:norm_B_bound}),
\begin{equation*}
|v_k^h|_{\alpha + \bar{\beta}} \le C \big| f_k^h - \mathcal{A}^{\alpha, k, h} v - \mathcal{B}^{\alpha, k, h} v\big|_{\bar{\beta}} \le C |f|_\beta, k = 1, \ldots, d,
\end{equation*}
with $C$ being independent of $h$ and $f$. Also by Corollary~\ref{cor:operator_A_bound} and Lemma~\ref{lem:operator_B_bound},
\begin{equation} \label{eqn:norm_A & B_bound}
|\mathcal{B}^\alpha v_k^h|_{\bar{\beta}} \le C |v_k^h|_{\alpha + \bar{\beta}} \le C |f|_\beta \quad \mbox{and} \quad |\mathcal{A}^\alpha v_k^h|_{\bar{\beta}} \le C |v_k^h|_{\alpha + \bar{\beta}} \le C |f|_\beta.
\end{equation}

Therefore by (\ref{eqn:cauchy_diff}), for any $(t, x) \in \mathds{H}$,
\begin{eqnarray*}
v_k^h(t, x) - v_k^h(s, x) & = & \int_s^t \big[f_k^h(r, x) - \mathcal{A}_x^{\alpha, k, h} v(r, x + h e_k) - \mathcal{B}_x^{\alpha, k, h} v(r, x + h e_k)\big] dr \\
 & & - \int_s^t \big(\mathcal{A}_x^\alpha + \mathcal{B}_x^\alpha\big) v_k^h(r, x)dr, 0 \le s < t \le T, k = 1, \ldots, d
\end{eqnarray*}
and by (\ref{eqn:norm_f_bound}) - (\ref{eqn:norm_A & B_bound}),
\begin{eqnarray*}
|v_k^h(t, x) - v_k^h(s, x)| & \le & \big(\big| f_k^h - \mathcal{A}_x^{\alpha, k, h} v - \mathcal{B}_x^{\alpha, k, h} v\big|_{\bar{\beta}} + \big| \big(\mathcal{A}_x^\alpha + \mathcal{B}_x^\alpha\big) v_k^h\big|_{\bar{\beta}}\big)|t - s| \\
 & \le & C |t - s|,
\end{eqnarray*}
where $C$ is a constant independent of $h$ and $k$. Hence, $v_k^h(t, x), k = 1, \ldots, d$ are equicontinuous in $(t, x)$ and by the Arzel\`{a}-Ascoli theorem, for each $h_n \to 0$, there exist a subsequence $\{h_{n_j}\}$ and continuous functions $v_k(t, x), (t, x) \in \mathds{H}, k = 1, \ldots, d$, such that $v_k^{h_{n_j}}(t, x) \to v_k(t, x)$ uniformly on compact subsets of $\mathds{H}$ as $j \to \infty$. By Lemma~\ref{lem:P2Limit}, $v_k \in C^{\alpha + \bar{\beta}}$ and $|v_k|_{\alpha + \bar{\beta}} \le C |f|_\beta, k = 1, \ldots, d$.

It then follows from passing to the limit in (\ref{eqn:cauchy_diff}) and the dominated convergence theorem that $u_k$ is the unique solution to
\begin{eqnarray*}
\big(\partial_t + \mathcal{A}_x^\alpha + \mathcal{B}_x^\alpha\big) v_k(t, x) & = & \partial_k f(t, x) - \big(\partial_k \mathcal{A}_x^\alpha\big) v(t, x) - \big(\partial_k \mathcal{B}_x^\alpha\big) v(t, x), \\
v_k(T, x) & = & 0, k = 1, \dots, d
\end{eqnarray*}
and so $v_k^{h_n}(t, x) \to v_k(t, x), \forall h_n \to 0$. Hence,
\begin{equation*}
v_k(t, x) = \lim_{h \to 0} v_k^h(t, x) = \lim_{h \to 0} \frac{v(t, x + h e_k) - v(t, x)}{h} = \partial_k v(t, x),
\end{equation*}
$\partial_k v \in C^{\alpha + \bar{\beta}}(\mathds{H}), k = 1, \dots, d$, and $|\nabla v|_{\alpha + \bar{\beta}} \le C |f|_\beta $. Therefore, $v \in C^{\alpha + \beta}(\mathds{H})$ and the statement follows.
\end{proof}

With Theorem~\ref{thm:solution_Kolmogorov_diffusion_jump}, Corollary~\ref{cor:solution_Cauchy_diffusion_jump} follows straightforward and the proof is omitted here.

\begin{remark} \label{rem:norm_partial_of_time_diffusion_jump}
If the assumptions of Corollary~\ref{cor:solution_Cauchy_diffusion_jump} hold and $v \in C^{\alpha + \beta}(\mathds{H})$ is the solution to $(\ref{eqn:Cauchy_diffusion_jump})$, then $\partial_t v = f - \mathcal{A}_x^\alpha v - \mathcal{B}_x^\alpha v$ and by Corollary~\ref{cor:operator_A_bound} and Lemma~\ref{lem:operator_B_bound}, $|\partial_t v|_\beta \le C (|f|_\beta + |g|_{\alpha + \beta})$.
\end{remark}

\subsection{One-Step Estimate}

To determine the rate of convergence of the Euler approximation, a key step is to estimate the conditional expectation of each increment. Lemma~\ref{lem:one_step_estimate_diffusion_jump} provides such a one-step estimate.

\begin{lemma} \label{lem:one_step_estimate_diffusion_jump}
Let $Y = \{Y_t\}_{t \in (0, T]}$ be the weak Euler approximation with step size $\delta \in (0, 1)$ of the stochastic process $X = \{X_t\}_{t \in (0, T]}$ defined by $(\ref{eqn:diffusion_jump})$. For $\alpha \in (0, 2]$ and $\beta > 0, \beta \ne \alpha, \beta \notin \mathds{N}$, assume \textup{A1} and \textup{A2} hold. Then there exists a constant $C$ such that for all $f \in C^\beta(\mathds{R}^d)$,
\begin{equation*}
\big|\mathrm{E}\big[f(Y_s) - f(Y_{\tau_{i_s}}) | \mathcal{F}_{\tau_{i_s}} \big] \big| \le C |f|_\beta \delta^{\kappa(\alpha, \beta)}, \forall s \in [0, T],
\end{equation*}
where $i_s = i$ if $\tau_i \le s < \tau_{i + 1}$ and $\kappa(\alpha, \beta)$ is as defined in Theorem~\ref{thm:main_diffusion_jump}.
\end{lemma}

The proof of Lemma~\ref{lem:one_step_estimate_diffusion_jump} is based on applying Ito's formula to $f(Y_s) - f(Y_{\tau_{i_s}})$, $f \in C^\beta(\mathds{R}^d)$ and using the estimates for $\mathcal{A}_z f^\varepsilon$ and $\mathcal{B}_z f^\varepsilon$, which are stated in Lemma~\ref{lem:estimate_main_part_diffusion_jump} and Corollary~\ref{cor:estimate_lower_part_diffusion_jump}, respectively. Here for $\varepsilon \in (0, 1)$, $f^\varepsilon$ is the convolution defined by
\begin{equation*}
f^\varepsilon(x) = \int f(y) w^\varepsilon(x - y)dy = \int f(x - y) w^\varepsilon(y)dy, x \in \mathds{R}^d,
\end{equation*}
where $w^\varepsilon(x) = \varepsilon^{-d} w\big(\frac{x}{\varepsilon}\big)$, given a nonnegative smooth function $w \in C_0^\infty(\mathds{R}^d)$ with support on $\{|x| \le 1\}$ such that $w(x) = w(|x|)$, $x \in \mathds{R}^d$ and $\int w(x)dx = 1$.

Similar steps as those in \cite{MiZ11} can be followed to prove Lemma~\ref{lem:one_step_estimate_diffusion_jump}, as well as Lemma~\ref{lem:estimate_main_part_diffusion_jump} and Corollary~\ref{cor:estimate_lower_part_diffusion_jump}. The details are omitted here.

\begin{lemma} \label{lem:estimate_main_part_diffusion_jump}
Let $\alpha \in (0, 2)$, $\beta < \alpha$, $\beta \ne 1$, $\varepsilon \in (0, 1)$, and $f \in C^\beta(\mathds{R}^d)$. Then
\begin{enumerate}
\item[\textup{(i)}] there exists a constant $C$ such that
\begin{equation*}
|f^\varepsilon(x) - f(x)| \le C |f|_\beta \varepsilon^\beta, \forall x \in \mathds{R}^d;
\end{equation*}
\item[\textup{(ii)}] there exists a constant $C$ such that
\begin{equation*}
|\mathcal{A}_z^\alpha f^\varepsilon(x)| \le C |f|_\beta \varepsilon^{-\alpha + \beta}, \forall z, x \in \mathds{R}^d
\end{equation*}
and in particular,
\begin{equation*}
|\partial^\alpha f^\varepsilon(x)| \le C |f|_\beta \varepsilon^{-\alpha + \beta}, \forall x \in \mathds{R}^d;
\end{equation*}
\item[\textup{(iii)}] there exist constants $C$s such that for $k, l = 1, \ldots, d$,
\begin{eqnarray*}
|\partial_k f^\varepsilon(x)| & \le & C |f|_\beta \varepsilon^{-1 + \beta}, \forall x \in \mathds{R}^d, \mbox{ if } \beta < 1, \notag \\
|f^\varepsilon|_1 & \le & C |f|_1, \\
|\partial_{kl}^2 f^\varepsilon(x)| & \le & C |f|_\beta \varepsilon^{-2 + \beta}, \forall x \in \mathds{R}^d, \mbox{ if } \beta < 2, \notag
\end{eqnarray*}
and for $\alpha \in (1, 2)$,
\begin{eqnarray*}
|f^\varepsilon|_\alpha & \le & C |f|_\beta \varepsilon^{-\alpha + \beta}, \mbox{ if } \beta \in (0, 1], \notag \\
|\partial^{\alpha - 1} \nabla f^\varepsilon(x)| & \le & C |f|_\beta \varepsilon^{-\alpha + \beta}, \forall x \in \mathds{R}^d, \mbox{ if } \beta \in (1, \alpha).
\end{eqnarray*}
\end{enumerate}
\end{lemma}

\begin{corollary} \label{cor:estimate_lower_part_diffusion_jump}
Let $\varepsilon \in (0, 1)$. Assume $a(x)$ and $\displaystyle \int_{\mathds{R}_0^d} \big(|y|^\alpha \wedge 1\big) \rho^\alpha(x, y) \nu^\alpha(dy)$, $x \in \mathds{R}^d$ are bounded.
Then there exists a constant $C$ such that
\begin{equation*}
|\mathcal{B}_z^\alpha f^\varepsilon(x)| \le C |f|_\beta \varepsilon^{-\alpha + \beta}, \forall z, x \in \mathds{R}^d, \forall f \in C^\beta(\mathds{R}^d).
\end{equation*}
\end{corollary}

\subsection{Rate of Convergence}

With results on the backward Kolmogorov equations and one-step estimates, the rate of convergence stated in Theorem~\ref{thm:main_diffusion_jump} is proved by applying It\^{o}'s formula. The idea is the same as that in \cite{MiZ11}. The main steps are provided here for reference.

\begin{proof}
Let $v \in C^{\alpha + \beta}(\mathds{H})$ be the unique solution to (\ref{eqn:Cauchy_diffusion_jump}) with $f = 0$. By It\^{o}'s formula and Remark~\ref{rmk:martingale_diffusion_jump},
\begin{equation*} \label{eqn:expect_terminal}
\mathrm{E}[v(0, X_0)] = \mathrm{E}[v(T, X_T)] = \mathrm{E}[g(X_T)].
\end{equation*}

By Corollaries~\ref{cor:solution_Cauchy_diffusion_jump} and~\ref{cor:operator_A_bound}, Lemma~\ref{lem:operator_B_bound}, and Remark~\ref{rem:norm_partial_of_time_diffusion_jump}, for $s \in [0, T]$,
\begin{eqnarray*}
|\mathcal{A}_z^\alpha v(s, \cdot)|_\beta + |\mathcal{B}_z^\alpha v(s, \cdot)|_\beta & \le & C |v|_{\alpha + \beta} \le C |g|_{\alpha + \beta}, \\
|\partial_t v(s, \cdot)|_\beta & \le & C |g|_{\alpha + \beta}. \notag
\end{eqnarray*}

Then, by It\^{o}'s formula, Remark~\ref{rmk:martingale_diffusion_jump}, and Corollary~\ref{cor:solution_Cauchy_diffusion_jump}, it follows that
\begin{eqnarray*}
\mathrm{E}[g(Y_T)] - \mathrm{E}[g(X_T)] & = & \mathrm{E}[v(T, Y_T)] - \mathrm{E}[v(0, X_0)] \\
 & = & \mathrm{E}[v(T, Y_T) - v(0, Y_0)] \\
 & = & \mathrm{E}\big[\int_0^T [\partial_t v(s, Y_s) + \mathcal{A}_{Y_{\tau_{i_s}}}^\alpha v(s, Y_s) + \mathcal{B}_{Y_{\tau_{i_s}}}^\alpha v(s, Y_s)]ds\big] \\
 & = & \mathrm{E}\big[\int_0^T \big\{[\partial_t v(s, Y_s) - \partial_t v(s, Y_{\tau_{i_s}})] \\
 & & + [\mathcal{A}_{Y_{\tau_{i_s}}}^\alpha v(s, Y_s) - \mathcal{A}_{Y_{\tau_{i_s}}}^\alpha v(s, Y_{\tau_{i_s}})] \\
 & & + [\mathcal{B}_{Y_{\tau_{i_s}}}^\alpha v(s, Y_s) - \mathcal{B}_{Y_{\tau_{i_s}}}^\alpha v(s, Y_{\tau_{i_s}})]\big\}ds\big].
\end{eqnarray*}

Hence, by Lemma~\ref{lem:one_step_estimate_diffusion_jump}, there exists a constant $C$ independent of $g$ such that
\begin{equation*}
|\mathrm{E}[g(Y_T)] - \mathrm{E}[g(X_T)]| \le C |g|_{\alpha + \beta} \delta^{\kappa(\alpha, \beta)}.
\end{equation*}
The statement thus follows.
\end{proof}

\section{Conclusion}

The paper studies the weak Euler approximation of It\^{o} diffusion and jump processes. In particular, it investigates the dependence of convergence rate on the regularity of the coefficients and driving processes, under the assumption of H\"{o}lder continuity. The rate of convergence is estimated by first solving the associated backward Kolmogorov equation in H\"{o}lder space and then applying It\^{o}'s formula. It is proved that the rate of convergence is $\frac{\beta}{\alpha} \wedge 1$, given that the coefficients are $\beta$-H\"{o}lder continuous and the principal part of the jump intensity measure has a nondegenerate density with respect to the L\'{e}vy measure of a spherically-symmetric $\alpha$-stable process.

For the stochastic processes considered in this paper, a minor restriction on the scale of regularity is that $\beta \ne \alpha$ and $\beta \notin \mathds{N}$. A further step is to remove the restriction and to derive the rate of convergence for the whole H\"{o}lder-Zygmund scale.

\bibliographystyle{plain}
\bibliography{SDE_Refs}

\end{document}